\documentclass [12 pt, twoside]{amsart}
\usepackage{enumerate}

 \usepackage{amsmath,amsfonts,amsthm,amssymb}
\usepackage{color}

\newtheorem{thm}{Theorem}[section]
\newtheorem{lem}{Lemma}[section] 
\newtheorem{cor}{Corollary}[section]

\theoremstyle{definition}
\newtheorem*{Proof}{Proof}

\newcommand{\ld} {{\ldots}}
\newcommand{\sm} {{\smallsetminus}}

\newcommand{\de} {{\delta}}

\newcommand{\la} {{\lambda}}

\newcommand{\el} {{\ell}}
\newcommand{\e} {{\varepsilon}}

\newcommand{\dis}{\displaystyle}

\newcommand{\ch}{{\mathcal{H}}}

\newcommand{\ra}{{\rightarrow}}

\newcommand{\mfd}{{\mathfrak{D}}}

\newcommand{\qb}{$\quad\blacksquare$}
\def\1{\it1\hspace*{-0.150cm}{\footnotesize{I}}}

\def\R{{\mathbb{R}}}
\def\C{{\mathbb{C}}}

\def\N{{\mathbb{N}}}

\numberwithin{equation}{section}

\begin{document}

\title[ Non-existence of common hypercyclic functions ]{ Non-existence of common hypercyclic entire functions for certain families
of translation operators}

\author[George Costakis]{George Costakis}
\address{Department of Mathematics and Applied Mathematics, University of Crete, GR-700 13 Heraklion, Crete, Greece}
\email{costakis@uoc.gr}

\author[Nikos Tsirivas]{Nikos Tsirivas}
\address{Department of Mathematics and Applied Mathematics, University of Crete, GR-700 13 Heraklion, Crete, Greece}
\email{tsirivas@uoc.gr}

\author[Vagia Vlachou]{Vagia Vlachou}
\address{Department of Mathematics, University of Patras, 26500 Patras, Greece}
\email{vvlachou@math.upatras.gr}
\thanks{The second author during this research was supported by Science Foundation Ireland under Grant
09/RFP/MTH2149 }


\date{}


\begin{abstract}
Let $\ch(\C)$ be the set of entire functions endowed with the topology of local uniform convergence. Fix a sequence of non-zero complex numbers $(\la_n)$ with $\liminf_{n}\frac{|\la_{n+1}|}{|\la_n|}>2$. We prove that there exists no entire function $f$ such that for every $b\in \mathbb{C}\setminus \{ 0\}$ the set $\{ f(z+\la_nb): n=1,2,\ldots \}$ is dense in $\ch(\C)$. This, on one hand gives a negative answer to Question 2 in \cite{7} and on the other hand shows that certain results from \cite{Tsi}, \cite{Tsi2} are sharp.
\end{abstract}

\maketitle

\section{Introduction}

Let $\ch(\C)$ be the set of entire functions endowed with the topology of local uniform convergence. For every
$b\in\C\sm\{0\}$, let $T_b:\ch(\C)\ra\ch(\C)$ be the translation operator defined by the formula
$T_b(f)(z)=f(z+b)$ for $f\in\ch(\C)$, $z\in\C$. An old result due to Birkhoff \cite{4} 
says, in non-technical terms, that the operator $T_1$, although linear, exhibits quite complicated dynamical behavior. To be more
precise, there exists an entire function $f$ so that its iterates under $T_1$, $\{ T_1^n(f):n=1,2,\ldots \}$
(the symbol $T_1^n$ is understood as composing $T_1$ with itself $n$-times), form a dense set in $\ch(\C)$. This
means that the set $\{ f(z+n):n=1,2,\ldots \}$ is dense in $\ch(\C)$.

In modern terminology, this type of behavior is a particular instance of a general phenomenon in (functional)
analysis, so called \textit{hypercyclicity}, which appears frequently in any ``reasonable" topological vector
space. To recall, briefly, let $X$ be a real or complex separable topological vector space. A sequence $(S_n)$
of linear and continuous operators on $X$ is called \textit{hypercyclic} provided there exists $x\in X$ so that
the set $\{ S_n(x):n=1,2,\ldots \}$ is dense in $X$. Then $x$ is called \textit{hypercyclic} for $(S_n)$ and
$HC(\{ S_n\})$ denotes the set of all hypercyclic vectors for $(S_n)$. If in the previous definition the
sequence $(S_n)$ comes from the iterates of a single operator $S$ then we simply say that $S$ is \textit{hypercyclic}, $x$ is \textit{hypercyclic} for $S$ and $HC(S)$ are the hypercyclic
vectors for $S$. For a detailed study on this subject we refer to the recent books \cite{2}, \cite{10}.

In this work we focus on translation operators on $\ch(\C)$ and in particular we contribute to the problem of
\textit{common} hypercyclic functions for such operators. Let us fix a sequence of non-zero complex numbers
$(\la_n)$ with $|\la_n|\to +\infty$. It is well known that for every $b\in \C\setminus \{ 0\}$ the sequence
$(T_{\la_nb} )$ is hypercyclic in $\ch(\C)$ and the set $HC(\{ T_{\la_nb} \})$ is $G_{\de}$, i.e. countable
intersection of open sets, and dense in $\ch(\C)$, see for instance \cite{Grosse1}, \cite{Grosse2} (of course,
this is an extension of Birkhoff's result mentioned above). Therefore, an appeal of Baire's category theorem
shows that the intersection $\bigcap_{b\in B}HC(\{ T_{\la_nb} \})$ is always non-empty, in fact $G_{\de}$ and
dense in $\ch(\C)$, whenever $B$ is a countable subset of $\C\sm \{ 0\}$. However, whether the above
intersection remains non-empty whenever $B$ is uncountable is already a non-trivial problem! Indeed, our point
of departure is the following question, which was raised by the first author in \cite{7}.



\smallskip
{\bf Question } \cite{7}. \textit{Let $(\la_n)$ be a sequence of non-zero complex numbers such that
$\dis\lim_{n\ra+\infty}|\la_n|=+\infty$. Is it true that} $\bigcap\limits_{b\in\C\sm\{0\}}HC(\{T_{\la_nb}\})\neq
\emptyset$? \smallskip

In many cases the above Question admits a positive answer, as for instance in the cases: $\la_n=n$ \cite{6},
$\la_n =n\log (n+1)$ (\cite{7}, \cite{Tsi2}), $\la_n=n^2$, $\la_n=n^3$ etc. \cite{Tsi}. The main result of this work is the
following
\begin{thm}\label{thm2.1}
Let $(\la_n)$ be a sequence of non-zero complex numbers such that $$\liminf_{n\to +\infty}\frac{|\la_{n+1}|}{|\la_n|}>2.$$
Then for every non-degenerate line segment $I\subseteq\C\sm\{0\}$ the set $\bigcap\limits_{b\in I}HC(\{T_{\la_nb}\})$ is empty.
In particular, $\bigcap\limits_{b\in\C\sm\{0\}}HC(\{T_{\la_nb}\})=\emptyset$.
\end{thm}
Theorem \ref{thm2.1} has two main consequences. On one hand, it gives a negative answer to the above Question.
On the other hand, it sharpens certain recent results of the second author; for instance, the main result in
\cite{Tsi} implies the following:
$$
\bigcap_{b\in \mathbb{C}\setminus \{ 0\} }HC(\{ T_{e^{n^c}b} \}) \neq \emptyset \ \ \text{for every} \ \ 0<c<1.
$$
Now, we stress that the allowed growth $e^{n^c}$, $0<c<1$ in the previously mentioned result is optimal, as far
as the highest power $c$ concerns us, since
$$
\bigcap_{b\in \mathbb{C}\setminus \{ 0\} }HC(\{ T_{e^nb} \}) =\emptyset ,
$$
by Theorem \ref{thm2.1}. Observe that in order to conclude the lack of common hypercyclic functions for families
of translation operators as above, it is necessary to impose some kind of geometrical growth on the sequence
$(\la_n)$ in view of the following: for a sequence of non-zero complex numbers $(\la_n)$ with $|\la_n|\to
+\infty$,
$$\text{if}  \ \lim_{n\to
+\infty}\frac{|\la_{n+1}|}{|\la_n |}=1  \ \text{then}  \ \bigcap\limits_{b\in \C\sm\{0\}}HC(\{T_{\la_nb}\}) \neq
\emptyset ,$$ which is a particular case of the main result in \cite{Tsi}. To complete the picture, it remains
to deal with the case $1<\liminf_{n}\frac{|\la_{n+1}|}{|\la_n|}\leq 2$, which for us, at least up to now, is a
``grey zone" and perhaps reflects the limitations of our method.


One can find in the literature several results supporting the existence of common hypercyclic vectors for many
kinds of operators, besides translations, such as weighted shifts, adjoints of multiplication operators,
differentiation and composition operators; see for instance, \cite{AbGo}-\cite{Ber},
\cite{ChaSa1}-\cite{GaPa}, \cite{10}-\cite{13}. For results showing the
lack of common hypercyclic vectors for certain families of hypercyclic operators, we refer to \cite{Ba3}, \cite{BaGr1}, \cite{13}.

\section{Some combinatorial lemmas for intervals}
For the proof of Theorem \ref{thm2.1} we establish a series of lemmas which, we believe, are of independent interest. Let us introduce some standard notation. For an interval $I$, either in $\R$ or in $\C$, we denote its length by $|I|$. Let $V$ be a subset of $\R$ or $\C$. The symbol $\textrm{diam}V$ stands for the diameter of $V$.    

\subsection{Intervals in $\R$}

\begin{lem} \label{L1}
Let $I\subseteq\R$ be a compact interval and $V\subseteq I$ be open (in the usual topology of $\R$). Suppose, in addition, that there exist two positive numbers $a, A$ such that $2a\leq A$ and for every $x,\widetilde{x}\in V$ either $|x-\widetilde{x}|<a$ or $|x-\widetilde{x}|>A$. Then at least one of the following is true :
\begin{itemize}
  \item[(i)] $\textrm{diam}V \leq a$. 
  \item[(ii)] $I\sm V$ contains a closed interval of length $A$.
\end{itemize}
\end{lem}
\begin{Proof}
Consider the relation $\mfd\subseteqq V\times V$ defined by the following rule:
\[
(x,\widetilde{x})\in\mfd \ \ \text{if and only if} \ \ |x-\widetilde{x}|<a.
\]
Obviously, the relation $\mfd$ is reflexive and symmetric. We shall prove that it is also transitive. Let
$x_1,x_2,x_3\in V$ such that $(x_1,x_2)\in \mfd$ and $(x_2,x_3)\in \mfd$. Then $|x_1-x_2|<a$, $|x_2-x_3|<a$
and
\[
|x_1-x_3|\le|x_1-x_2|+|x_2-x_3|<2a\le A.
\]
Since $x_1,x_3\in V$, our hypothesis implies that either $|x_1-x_3|<a$ or $|x_1-x_3|>A$. In view of the above we conclude that  $|x_1-x_3|<a$, which shows that the relation $\mfd$ is transitive. Thus, $\mfd$ is an equivalence relation. Therefore we may consider  $V_j$, $j\in G$ the classes of equivalence, where $G$ is a set of indices. \smallskip

Then the following hold:
\begin{enumerate}
\item[(i)] The set $V_i$ is open in $\mathbb{R}$, for every $i\in G$.
\item[(ii)] $dist(V_i,V_j)\geq A$ for every $i,j\in G$ with $i\neq j$.
\item[(iii)] The set $G$ is finite.
\end{enumerate}

To prove (i) we fix $i\in G$ and it remains to prove that the set $V_{i}$ is open in $\R$. Let $x\in V_{i}$. Clearly, $x\in V$. Since $V$ is open in $\R$, there exists $\e>0$ such that $(x-\e,x+\e)\subseteq V$. Because the interval $(x-\e,x+\e)$ is connected and the function $f(y)=|x-y|$, $y\in (x-\e,x+\e )$ is continuous, the image of $f$ is connected. Thus, $\e \leq \frac{a}{2}$. Therefore $(x-\e,x+\e)$ is actually contained in $V_i$. \smallskip

We proceed with the proof of item (ii). Obviously, for every $x\in V_i$ and $y\in V_j$ we have  $|x-y|>A$ and the result follows.
\smallskip 

We are left with the proof of item (iii). Arguing by contradiction, suppose that the set $G$ is denumerable (we can not have an uncountable family of disjoint open sets in $\R$). Let $V_1, V_2, \ldots $ be the classes of equivalence and for every positive integer $n$ choose $x_n\in V_n$. Since $V_n\subset I$ for every $n=1,2,\ldots $ and $I$ is compact the sequence $(x_n)$ has a limit point. The last gives that $|x_{k_1}-x_{k_2}|<a$, for some $k_1, k_2\in \N$ with $k_1\neq k_2$ and so $x_{k_1}$, $x_{k_2}$ belong to the same equivalence class. This is clearly a contradiction. Therefore $G$ is a finite set. \smallskip

\textbf{Case I:} $G$ contains only one element. Then obviously $\textrm{diam} V\leq a$.

\textbf{Case II:} $G$ contains more than one element. Let us assume that $V_1,\ldots,V_n$ are the classes of equivalence, for some positive integer $n\geq 2$ .\\
Without loss of generality (changing the enumeration if necessary), we may assume that  $\displaystyle \inf V_{1}=\min_{i=1,2,\dots,n}\inf V_{i}$ and $\displaystyle \inf V_{2}=\min_{i=2,\dots,n}\inf V_{i}$ .\\
Then the closed interval $[\sup V_{1}, \sup V_{1}+A]$ is contained in $I\sm V$. Let us see why: 

\smallskip\noindent
$\bullet$  Since $V_1, V_2$ are  open, $\sup V_{1}\notin V_1$ and $\inf V_2 \notin V_2$.\\
$\bullet $ Since $I$ is closed,  $\sup V_{1},  \inf V_{2} \in I$.\\
$\bullet$ $|\inf V_2-\sup V_1|\geq A$ (because of (ii)).\\
$\bullet$  $\inf V_2-\sup V_{1} = (\inf V_2-\inf V_1)-(\sup  V_1-\inf V_1)\geq  A-a \geq 2a-a= a >0$.\\
$\bullet$  Since $I$ is an interval,  $[\sup  V_{1},  \inf  V_{2}]\subset I$.\\
$\bullet$ It is easy to see that $V$ does not intersect the  interval $[\sup V_{1},  \inf V_{2}]$.\\
The proof is complete  in view of $ [\sup V_{1}, \sup V_{1}+A]\subset [\sup V_{1},  \inf V_{2}] $. \qb 
\end{Proof}
The above lemma does not hold if we remove the assumption $A\geq 2a$. Let us give a counterexample. Assume that $a\leq A<2a$, then set $\e =\frac{2a-A}{3} >0$ and consider $I=[0,2a]$ and $V=(0,\e )\cup (a-\e ,a)\cup (2a-2\e ,2a-\e )$. If $A<a$ every open set satisfies the assumptions of the lemma, but not necessarily the conclusion.\medskip

The following corollary is actually the result we need for the proof of our main result. 
\begin{cor}\label{C}
Let $a,A$ be two positive numbers such that $2a\leq A$. Let, in addition, $I$ be a compact interval in $\R$ and $V\subset I$ be an open set. Assume that the following conditions hold:
\begin{itemize}
\item[(i)] $|I|\geq 2A+a$.
\item[(ii)] For every $x,y\in V$, either $|x-y|<a$ or $|x-y|>A$.
\end{itemize}
Then the set $I\setminus V$ contains a closed interval of length $A$.
\end{cor}
\begin{Proof}
In view of the previous lemma, we only have to deal with the case $\textrm{diam} V\leq a$. It suffices to prove that $I\setminus (\inf V,\sup V)$ has at least one connected component of length at least $A$ (note that this set has at most $2$ connected components). If this is not the case, then $|I|<2A+\sup V-\inf V$ and we have a contradiction.\qb
\end{Proof}


\subsection{Intervals in $\C $}
The next corollary is nothing but the complex version of Corollary \ref{C}.
\begin{cor}  \label{C2}
Let $a,A$ be two positive numbers, such that $2a\leq A$. Let, in addition, $I\subset \C$ be a line segment and $V\subset \C$ be an open set (in the usual topology of $\C$). Assume that the following conditions hold:
\begin{itemize}
\item[(i)] $|I|\geq 2A+a$.
\item[(ii)] For every $z,w\in V$, either $|z-w|<a$ or $|z-w|>A$.
\end{itemize}
Then the set $I\setminus V$ contains a line segment of length $A$.
\end{cor}
\begin{Proof}
\textbf{Case I:} If $I$ is contained in $\R$ and $I\cap V$ is open in $\R$ then the result follows from Corollary \ref{C}.

\textbf{Case II:} If $I$ is contained in $\R$ and $I\cap V$ is not open in $\R$, set $V^*=(I\cap V)\setminus \{ \textrm{endpoints} \,\,\, \textrm{of}\,\,\, I\} $ and apply Corollary \ref{C} to $V^*$. The result follows.

\textbf{Case III:} In the general case where $I \subset \C$ but $I$ is not contained in $\R$, we use an isometry $\phi :\C \ra \C$, composing a translation with a rotation so that $I':=\phi (I)$ is a closed interval in $\R$ and $|I'|=|I|$. By this, we transfer our problem to the previous cases and everything works nicely, since the isometry $\phi$ preserves closed, open sets and lengths. \qb 
\end{Proof}

\section{Two elementary lemmas}
We prove two elementary lemmas, which will connect the above results with our main result.

\smallskip

Let $k$ be a positive number, $\lambda \in  \C\sm\{0\}$ and $f:\C\ra\C$, $g: \{ z: |z|\leq k\}\ra\C$ be complex functions. We will use the following notations:
$$V_{\lambda }(f,g,k):=\left\{ b\in \C :\sup_{|z|\leq k}|f(z+\lambda b)-g(z)|<\frac{1}{k} \right\} ,$$
$$V_{\lambda }(f,k):=\left\{ b\in \C :\sup_{|z|\leq k}|f(z+\lambda b)-z|<\frac{1}{k} \right\} .$$

For $z\in \C$ and $r>0$, $D(z,r)$, $\overline{D(z,r)}$ denote the open and closed disk with center $z$ and radius $r$, respectively.  

\begin{lem} \label{L6}
Let $k$ be a positive number, $\la \in\C\sm\{0\}$ and $f:\C\ra\C$ a complex function. If
$b_1,b_2\in V_{\la}(f,k)$ we have:
\[
\text{either} \ \ |b_1-b_2|<\frac{2}{k|\la |} \ \ \text{or} \ \ |b_1-b_2|>\frac{2k}{|\la |}.
\]
\end{lem}

\begin{Proof}
Let $b_1,b_2\in V_{\la}(f,k)$.\smallskip

{\bf Case 1}. $\overline{D(\la b_1,k)}\cap\overline{D(\la b_2,k)}\neq\emptyset$.
\smallskip

Then there exist complex numbers $z_1$, $z_2$ with $|z_1|, |z_2|\leq k$ such that
$w:=z_1+\la b_1=z_2+\la b_2$.
Since $|f(w)-z_1|<1/k$, $|f(w)-z_2|<1/k$ we have $|z_1-z_2|<2/k$. Therefore, $|b_1-b_2|<2/(k|\la |)$. \smallskip

{\bf Case 2}. $\overline{D(\la b_1,k)}\cap\overline{D(\la b_2,k)}=\emptyset$. \smallskip

Then $|\la b_1-\la b_2|>2k$. This completes the proof of the lemma. \qb
\end{Proof}
\begin{lem} \label{L7}
Let $k$ be a positive number, $\la \in \C$ and $f:\C\ra\C$ be a continuous function. Fix any function $g:\{ z: |z|\leq k\}\ra\C$. Then the set $V_{\la}(f,g,k)$ is an open subset of $\C$.
\end{lem}
\begin{Proof}
Let $b\in V_{\la}(f,g,k)$. Since $f$ is continuous, there exists $\delta$ with $0<\delta <1$ such that whenever $z,w\in \overline{D(\la b,k+1)}$ and $|z-w|<\delta$ the distance between $f(z)$ and $f(w)$ is less than $(1/k)-\sup_{|z|\leq k}|f(z+\la b)-g(z)|$. Using the last and the triangle inequality it is straightforward to check that the open disk $D(b,\frac{\delta}{|\la |})$ is contained in $V_{\la }(f,g,k)$. This completes the proof.\qb
\end{Proof}

\section{ Proof of Theorem \ref{thm2.1}}

We fix a compact, non-degenerate interval $I\subseteq\C\sm\{0\}$. Arguing by contradiction, assume that there exists a function
$f\in\bigcap\limits_{b\in I}HC(\{T_{\la_nb}\})$. Let some $\el>2$ with
\[
2<\el<\underset{n\ra+\infty}{\lim\inf}\bigg|\frac{\la_{n+1}}{\la_n}\bigg|.
\]
Without loss of generality we may assume that $\Big|\dfrac{\la_{n+1}}{\la_n}\Big|>\el$ for every $n=1,2,\ldots
$. Fix a positive number $k$ satisfying:
\begin{equation} \label{eq1}
 k^2\ge 2,
\end{equation}

\begin{equation} \label{eq2}
 k^2\ge\dfrac{1}{\el-2}.
\end{equation}

Since $\dis\lim_{n\ra+\infty}|\la_n|=\infty$, there exists a positive integer $m$ such that
\begin{equation} \label{eq3}
\frac{4k}{|\la_n|}+\frac{2}{k|\la_n|}<|I|, \ \ \text{for every} \ \ n=m,m+1,\ldots .
\end{equation}
For every $n\in\N$ with $n\ge m$ we consider the set
\[
V_n:=V_{\la_n}(f,k)=\bigg\{b\in \mathbb{C}:\sup_{|z|\le k}|f(z+\la_nb)-z|<\frac{1}{k}\bigg\}.
\]
Since $f\in\bigcap\limits_{b\in I}HC(\{T_{\la_{n+m}b}\})$, keep in mind that $m$ is fixed,
the inclusion 
\begin{eqnarray*}
I\subset \bigcup^{+\infty}_{n=m}V_n
\end{eqnarray*}
is straightforward. Moreover, in view of Lemma \ref{L7}, the sets $V_n$, $n=1,2,\ldots $ are open in $\C$. Therefore, since $I$ is compact we have:
$$I=\bigcup\limits^N_{n=m}(V_n\cap I),$$
for some positive integer $N\ge m$. Let us define
\[
n_1:=\min\{n\in\N: m\le n\le N \ \ \text{and} \ \ I\cap V_n\neq\emptyset\}.
\]
By Lemma \ref{L6}, if $b_1, b_2\in V_{n_1}$ then either $|b_1-b_2|<2/(k|\la_{n_1}|)$ or $|b_1-b_2|>2k/|\la_{n_1}|$ and 
$$\frac{ \frac{2k}{|\la_{n_1}|} }{ \frac{2}{k|\la_{n_1}| }}=k^2\geq 2.$$
In view of (\ref{eq3}) we conclude that $I$, $V_{n_1}$ and $a:=\frac{2}{k|\la_{n_1}|}$, $A:=\frac{2k}{|\la_{n_1}|}$ satisfy the hypothesis of Corollary \ref{C2}. Hence, there exists a line segment $I_{n_1}$ such that: \smallskip

$I_{n_1}\subset I$,\smallskip

$I_{n_1}\cap V_{n_1}=\emptyset$, \smallskip

$|I_{n_1}|=A=\dfrac{2k}{|\la_{n_1}|}$. \smallskip

Thus, necessarily, 
$$I_{n_1}\subset\bigcup\limits^N_{n=n_1+1}V_n.$$ 
We iterate this argument. Let us define
\[
n_2:=\min\{n\in\N:n_1+1\le n\le N \ \ \text{and} \ \ I_{n_1}\cap V_n\neq\emptyset\}.
\]
Moreover, by (\ref{eq2}) we get
\begin{align*}
k^2\ge\frac{1}{\el-2}&\Leftrightarrow(\el-2)k\ge\frac{1}{k}\Leftrightarrow
\el k\ge 2k+\frac{1}{k} \\
&\Leftrightarrow\frac{2\el k}{|\la_{n_2}|}\ge\frac{4k}{|\la_{n_2}|}+\frac{2}{k|\la_{n_2}|}.
\end{align*}
Therefore,
\begin{equation} \label{eq4}
|I_{n_1}|=\bigg|\frac{2k}{\la_{n_1}}\bigg|>\frac{2\el
k}{|\la_{n_2}|}\ge\frac{4k}{|\la_{n_2}|}+\frac{2}{k|\la_{n_2}|}.
\end{equation}
Using (\ref{eq4}) and arguing as above we see that all the hypothesis of Corollary \ref{C2} are fulfilled for
$I:=I_{n_1}$, $V:=V_{n_2}$, $a:=\frac{2}{k|\la_{n_2}|}$ and $A:=\frac{2k}{|\la_{n_2}|}$. Hence, there exists a line
segment $I_{n_2}$, such that: \smallskip

$I_{n_2}\subset I_{n_1}$, \smallskip

$I_{n_2}\cap V_{n_2}=\emptyset$, \smallskip

$|I_{n_2}|=\dfrac{2k}{|\la_{n_2}|}$. \smallskip

The above imply that
$$
I_{n_2}\cap V_{n_2}=\emptyset \ \ \text{and} \ \ I_{n_2}\subset\bigcup^N_{n=n_2+1}V_n.
$$
Continuing this process, after a finite number of steps (at most $N-m+1$) we will end up with a subsegment of
$I$ disjoint from all $V_n$ for $n\in\{m,\ld,N\}$ and this is obviously a contradiction. \qb


\end{document}